\documentclass[12pt]{article}
\usepackage[scaled]{helvet} 
\usepackage{zlmtt} 
\usepackage[T1]{fontenc}
\usepackage{textcomp}
\normalfont
\usepackage{amssymb}
\usepackage{amsfonts}
\usepackage{amsmath,amssymb, amsthm}
\usepackage{euscript}
\usepackage{enumerate}
\usepackage{verbatim}
\usepackage{mathrsfs, color}
\usepackage{multicol}

\usepackage[latin1]{inputenc}
\usepackage[all]{xy}

\newtheorem{thm}{Theorem}[section]
\newtheorem{thm*}{Theorem }
\newtheorem{prop}[thm]{Proposition}
\newtheorem{lem}[thm]{Lemma}
\newtheorem{Cor}[thm]{Corollary}

\usepackage[all]{xy}

\def\bee{\begin{eqnarray}}
\def\bes{\begin{eqnarray*}}
\def\eee{\end{eqnarray}}
\def\ees{\end{eqnarray*}}
\def\a{\alpha}
\def\b{\beta}

\def\l{\lambda}

\def\s{\sigma}

\def\Proof{{\it Proof.\ }\ }
\def\ctd{\hfill$\Box$}

\def\VV{\mathcal V}
\def\A{\mathcal A}

\def\e{\varepsilon}
\def\Skew{\operatorname{Skew\,}}
\def\Alt{{Alt}}

\title{Skew-symmetric identities of finitely generated alternative  algebras}

\author{Ivan P. Shestakov\footnote{Supported by CNPq (Brazil), process 304313/2019-0, and by
FAPESP (Brazil), process 2018/23690-6}}

\date{}

\begin{document}
\maketitle

\abstract{We prove that for every natural number $n$  there exists a natural number $N(n)$ such that {\em every} multilinear skew-symmetric polynomial on $N(n)$ or more  variables which vanishes in the free associative algebra vanishes as well in any $n$-generated alternative algebra over a field of characteristic $0$. Before this was proved only for a series of such polynomials  constructed by the author in \cite{She-77}}. 

\section{Introduction}
\hspace{\parindent}
Let $\VV$ be a variety of algebras. Denote by $\VV_n$ the subvariety of $\VV$ generated by the $\VV$-free algebra on $n$ free generators; then we have 
\bes
\VV_1\subseteq \VV_2\subseteq\cdots\subseteq \VV_n\subseteq \cdots,\ \ \VV=\cup_{n}\VV_n.
\ees
A minimal number $n$ for which $\VV_n=\VV$ (if such a number exists) is called {\em the basic rank} of the variety $\VV$ and is denoted as $r_b(\VV)$ (see {\cite{Ma-70}). If $\VV\neq \VV_n$ for any $n$ then the basic rank of $\VV$ is called to be infinite. 

It is well known that for the varieties  of associative and Lie  algebras the basic rank is equal to 2. A.I.Shirshov  in \cite{Dniestr} posed a problem on the basic rank of the varieties of alternative $Alt$, Jordan $Jord$, Malcev $Malc$, and other algebras. 
In \cite{She-77} the author constructed a series of  multilinear skew-symmetric polynomials $g_m$ on $2m+1$ variables over a field of characteristic $\neq 2$ which vanish in any $n$-generated alternative or Malcev algebra for $m>[\tfrac{n^3+5n}{12}]$, but do not vanish in free alternative and free Malcev algebras of infinite rank. In particular, this implies that $r_b(Alt)$ and $r_b(Mal)$ are infinite.

Here we  prove more general result for alternative algebras. Namely, we prove that for any natural number $n$ there exists a number  $N=N(n)$, such that 
{\bf EVERY} multilinear skew-symmetric polynomial on $N(n)$ variables which vanishes in the free associative algebra vanishes in any $n$-generated alternative algebra. 

A similar result for Malcev algebras was proved in \cite{She-16}.

\section{Skew-symmetric functions in flexible algebras}
Let $F$ be a field of characteristic not 2.
Recall that an algebra $A$ over $F$ is called {\em flexible} if it satisfies the identity
\[
(x,y,x)=0
\]
where $(x,y,z)=(xy)z-x(yz)$ denotes the {\em associator} of elements $x,y,z$.
Any flexible algebra satisfies  the identity \cite[Lemma 1]{She-71} 
\bee\label{id_flex}
[x^2,y]=x\circ [x,y],
\eee
where $x\circ y=xy+yx$ denotes the {\em Jordan product} of $x,y$. 
We construct here a new series of skew-symmetric polynomials $f_m$ on $m$ variables in flexible algebras that have properties similar to those of polynomials $g_m$ constructed in \cite{She-77}.

Set $f_1(x)=x,\ f_2(x_1,x_2)=[x_1,x_2]$, and then, by induction, 
\[
f_{m+1}(x_1,x_2,\ldots,x_{m+1})=\sum_{1\leq i<j\leq m+1}(-1)^{i+j-1}f_m([x_i,x_j],x_1,\ldots,\check{x}_i,\ldots,\check{x}_j,\ldots,x_{m+1}),
\]
where $\check{x}_k$ means that the element $x_k$ is omitted.
\begin{lem}\label{lem_skew-symmetric}
Polynomials $f_m(x_1,\ldots,x_m)$ for $m>1$  are skew-symmetric on its variables, that is, they vanish if any two of arguments $x_i$ coincide. 
\end{lem}
\Proof
The statement is true for $m=2$. Assume now that it is true for certain $m\geq 2$ and prove it for $m+1$.
Let $x_i=x_j,\, i<j$. Then we have, by induction,
\bes
f_{m+1}(x_1,\ldots,x_{m+1})&=&\sum_{k=1}^{i-1}(-1)^{k+i-1}f_m([x_k,x_i],x_1,\ldots,\check{x}_k,\ldots,\check{x}_i,\ldots,x_{m+1})\\
&+&\sum_{j\neq k=i+1}^{m+1}(-1)^{i+k-1}f_m([x_i,x_k],x_1,\ldots,\check{x}_i,\ldots,\check{x}_k,\ldots,x_{m+1})\\
&+&\sum_{i\neq k=1}^{j-1}(-1)^{k+j-1}f_m([x_k,x_j],x_1,\ldots,\check{x}_k,\ldots,\check{x}_j,\ldots,x_{m+1})\\
&+&\sum_{k=j+1}^{m+1}(-1)^{j+k-1}f_m([x_j,x_k],x_1,\ldots,\check{x}_j,\ldots,\check{x}_k,\ldots,x_{m+1}).
\ees  
Now, if $k<i$, we have
\bes
f_m([x_k,x_i],x_1,\ldots,\check{x}_k,\ldots,\check{x}_i,\ldots,x_j,\ldots,x_{m+1})=\\
(-1)^{j-i-1}f_m([x_k,x_j],x_1,\ldots,\check{x}_k,\ldots,{x}_i,\ldots,\check{x}_j,\ldots,x_{m+1})
\ees
The similar equality we have for $k>j$. Finally, for $i<k<j$ we have
\bes
f_m([x_i,x_k],x_1,\ldots,\check{x}_i,\ldots,\check{x}_k,\ldots,x_j,\ldots,x_{m+1})=\\
(-1)^{j-i-2}f_m([x_k,x_j],x_1,\ldots,{x}_i,\ldots,\check{x}_k,\ldots,\check{x}_j,\ldots,x_{m+1})=\\
(-1)^{j-i-1}f_m([x_j,x_k],x_1,\ldots,{x}_i,\ldots,\check{x}_k,\ldots,\check{x}_j,\ldots,x_{m+1}).
\ees
Therefore, the first two sums in the expression for $f_{m+1}$ just differ by sign from the  last two sums, and the whole sum is 0.
\ctd
\begin{lem}\label{lem_sskew-symmetric} 
Polynomials $f_m$ are {strongly skew-symmetric}, that is, they satisfy the identity
\bee\label{id_sskew}
f_m(x^2,x,x_1,\ldots,x_{m-2})=0.
\eee
\end{lem} 
\Proof
The statement is true for $m=2$; assume that it is true for some $m\geq 2$ and prove for $m+1$.
Since the based filed $F$ has at least 3 elements,  identity (\ref{id_sskew}) implies 
\bee\label{id_sskew1}
f_m(x^2,y,x_1,\ldots,x_{m-2})+f_m(x\circ y,x,x_1,\ldots,x_{m-2})=0.
\eee
By the induction assumption and skew-symmetricity of $f_m$ we have
\bes
f_{m+1}(x^2,x,x_1,\ldots,x_{m-1})&=&\\
\sum_{i=1}^{m-1}(f([x^2,x_i],x,x_1,\ldots,\check{x}_i,\ldots,x_{m-1})&-&f([x,x_i],x^2,x_1,\ldots,\check{x}_i,\ldots,x_{m-1}) )
\ees
Denote, for shortness, $f(a,b)=f_{m}(a,b,x_1,\ldots,\check{x}_i,\ldots,x_{m-1})$. Then we have by (\ref{id_flex}) and (\ref{id_sskew1})
\bes
f([x^2,x_i],x)=f(x\circ[x,x_i],x)=-f(x^2,[x,x_i)=f([x,x_i],x^2),
\ees
which proves that $f_{m+1}(x^2,x,x_1,\ldots,x_{m-1})=0$.  

\ctd 

The following theorem essentially generalizes the corresponding result from \cite{She-77}.
\begin{thm}
Let $\VV$ be a variety of flexible algebras over a field of characteristic $\neq 2$ such that every finitely generated anticommutative algebra in $\VV$ is nilpotent.
Then for any $m$ there exists $N=N(m)$ such that any $m$-generated algebra in $\VV$ satisfies the identity $f_N=0$. 
If, moreover, for any $k$ polynomial $f_k$ is non-zero in the free $\VV$-algebra of countable rank then $r_b(\VV)$ is infinite.
The last condition holds if the characteristic of the based field is zero and $\VV$ contains the variety of alternative algebras. 
\end{thm}
\Proof
Let $A$ be the free $\VV$-algebra on $m$ generators. Denote by $I_2(A)$ the subspace of $A$ generated by all the squares of elements of $A$. By \cite{She-71}, $I_2(A)$ is an ideal of $A$ and the quotient algebra $\bar A=A/I_2(A)$ is anticommutative. Since $A$ is $m$-generated, so is $\bar A$; therefore, $\bar A$ is nilpotent.    Let $\bar A^n=0$, then $A^n\subseteq I_2(A)$, that is, every monomial $u$ from $A$ with $\deg(u)\geq n$ may be written as a linear combination of Jordan products $u_i\circ v_i$ where $\deg u_i+\deg v_i=\deg u$. 
Denote by $M_k$ a subspace of $A$ generated by the set $\{f_{k+1}(u_1,\ldots,u_k,u)\}$, where $u_i, u$ are monomials on generators and $\deg u_i<n$. Let us show that $f_{k+1}(A,\ldots,A)\subseteq M_k$. Linearizing identity \eqref{id_sskew1} we get in view of skew-symmetry of $f_{k+1}$
\bee\label{id_sskew2}
f_{k+1}(x\circ y,x_1,\ldots,x_{k-1},z)=f_{k+1}(x,x_1,\ldots,x_{k-1},z\circ y)+f_{k+1}(y,x_1,\ldots,x_{k-1},z\circ x).
\eee
Consider $f_{k+1}(a_1,\ldots,a_k,a)$ where $a_i,\,a$ are monomials on generators. If $\deg a_1\geq n$ then $a_1$ may be represented as  a linear combination of Jordan products $u\circ v$ where $\deg u,\deg v<\deg a_1$. By \eqref{id_sskew2}, then $f_{k+1}(a_1,\ldots,a_k,a)$ is a linear combination of elements $f_{k+1}(u,a_2,\ldots,a_k,a')$ with $\deg u<\deg a_1$. Continuing in this way, we eventually obtain that $f_{k+1}(a_1,\ldots,a_k,a)\in M_k$, which proves that $f_{k+1}(A,\ldots,A,A)\subseteq M_k$.

Let now $N=\dim A/I_2(A)$, prove that  $M_{N+1}=0$. In fact, since $A^n\subseteq I_2(A)$, any $N+1$ monomials $u_1,\ldots,u_{N+1}$ of degree less than $n$ in $A$ are linearly depending, hence in view of skew-symmetricity, $f_{N+2}(u_1,\ldots,u_{N+1},A)=0$ and $M_{N+1}=0$. By the previous considerations, this implies that $f_{N+2}(A,\ldots,A)=0$.

Therefore, if the elements $f_{k}(x_1,\ldots,x_k)\neq 0$ in the free $\VV$-algebra of infinite rank then $r_b(\VV)$ is infinite.

Finally, let $F$ be a field of characteristic zero. The similar proof as in \cite[Chapter 13]{ZSSS} shows that polynomials $f_k$ are nonzero in free metabelian alternative algebra on $k$ or more generators. Therefore, they are nonzero in the free $\VV$-algebra on $k$ or more generators.

\ctd

Recall that a flexible algebra is called {\em noncommutaive Jordan} if it satisfies the identity
\[
(x^2,y,x)=0,
\]

\begin{Cor}\label{cor1}
Let $\VV$ be a variety of noncommutative Jordan algebras over a field of characteristic zero defined by  the identity
\bes
([x,y],z,z)&=&0.
\ees
Then the basic rank $r_b(\VV)$ is infinite.
\end{Cor}
\Proof
In fact, it was proved in \cite{She-71} that  every finitely generated anticommutative algebra in  this variety is nilpotent.
It is also clear that this variety contains the variety of alternative algebras.

\ctd

\begin{Cor}\label{cor2}
Every $m$-generated  alternative algebra over a field of characteristic $\neq 2$ satisfies the identity $f_{n+1}=0$ for $n>1+m+\left(\begin{array}{c} m\\2\end{array}\right)+\left(\begin{array}{c} m\\3\end{array}\right)$.
\end{Cor}
\Proof
In fact, it is well known (see \cite[p.122]{ZSSS}) that an anticommutative alternative algebra $A$ over a field of characteristic $\neq 2$  with a set of generators $\{\a_i,\, i\in I\}$ is nilpotent of index 4 and is spanned by the elements of the type  $a_i,\,a_ia_j,\, (a_ia_j)a_k$ where $i<j<k$.  
Therefore, if $A$ has $m$ generators then  $\dim_F(A/I_2(A))\leq m+\left(\begin{array}{c} m\\2\end{array}\right)+\left(\begin{array}{c} m\\3\end{array}\right)$.   Hence we may take in this case $N(m)=m+\left(\begin{array}{c} m\\2\end{array}\right)+\left(\begin{array}{c} m\\3\end{array}\right)$.

\ctd

\section{Skew-symmetric elements and one-generated superalgebras}
\hspace{\parindent}
In this section we will assume that the field $F$ has characteristic $0$.

Let us recall the definition of a superalgebra in a given variety of algebras $\VV$.  In general, a superalgebra $A$ is just a ${\mathbf Z}_2$-graded algebra: $A=A_0\oplus A_1$, where $A_iA_j\subseteq A_{i+j\,(mod\,2)}$. Let $G={\rm alg}\langle 1,e_1,\ldots,e_n,\ldots|\, e_ie_j=-e_je_i\rangle $ be the Grassmann superalgebra, where $G_0$ is spanned by 1 and the even products $e_{i_1}e_{i_2}\cdots e_{i_k},\, 1<i_1<i_2<\cdots<i_k,\, k$ even, and $G_1$ is spanned by the odd products $e_{i_1}e_{i_2}\cdots e_{i_k},\, 1<i_1<i_2<\cdots<i_k,\, k$ odd. Then a superalgebra $A=A_0+A_1$ is called a {\em $\VV$-superalgebra} if its {\em Grassmann envelope}  $G(A)=G_0\otimes A_0+G_1\otimes A_1$, considered as an algebra,  belongs to $\VV$. 
When working with superalgebras, we always assume that considered elements are homogeneous, that is, {\em even} (belonging to $A_0$) or {\em odd} (belonging to $A_1$). 
If the defining identities of the variety $\VV$ are known, one can easily write the defining  {\em super-identites} for the variety of $\VV$-superalgebras. Recall that an algebra $A$ is called {\em alternative} (see {ZSSS}) if it satisfies the identities 
\bes
(x,y,y)&=&0,\\
(x,x,y)&=&0,
\ees
where $(x,y,z)=(xy)z-x(yz)$ denotes the {\em associator} of the elements $x,y,z$. Now, a superalgebra $A=A_0\oplus A_1$ is called an {\em alternative superalgebra} if it satisfies the superidentities
\bes
(x,y,z)+(-1)^{|y||z|}(x,z,y)&=&0, \\
(x,y,z)+(-1)^{|x||y|}(y,x,z)&=&0,
\ees
where $|z|=i$ if $z\in A_i$ denotes the parity of a homogeneous element $z$.
Denote by $[x, y]_s = xy +(-1)^{|x||y|}yx$ the {\em super-commutator} of the homogeneous
elements $x, y$, and by $x_s\circ y  = xy + (-1)^{|x||y|}yx$ their {\em super-Jordan product}.

Recall some results on relation between the space of skew-symmetric elements in a free $\VV$-algebra of countable rank and the free $\VV$-superalgebra $F_{\VV}[\emptyset;x]$ generated by one odd element $x$ (see {\cite{She-99, She-03, ShZh-06}).

Let $u(x)$ be a homogenious of degree $n$ element from $F_{\VV}[\emptyset;x]$. Write $u(x)$ in the form $u(x)=v(x,x,\ldots,x)$ where $v(x_1,\ldots,x_n)$ is a multilinear element. Set 
$$
{\rm  (Skew}\,u) (x_1,x_2,\ldots,x_n)=\sum_{\s\in \Sigma_n} (-1)^{sgn\, \s}v(x_{\s(1)},x_{\s(2)},\ldots,x_{\s(n)}).
$$
For example, if $u=xx$ then $({\rm Skew}\, u)(x_1,x_2)=[x_1,x_2]:=x_1x_2-x_2x_1$, the {\em  commutator} of elements $x_1,x_2$; if $u=(xx)x$ then $({\rm Skew}\, u)(x_1,x_2,x_3)=\sum_{\s\in\Sigma_3}(-1)^{sgn\,\s}(x_{\s(1)}x_{\s(2)})x_{\s(3)}$.

\begin{prop} \label{prop_u(x)}\cite{She-99, She-03, ShZh-06}.
Let $F$ be a field of characteristic $0$. The application $u\mapsto {\rm Skew}\, u$ establishes an isomorphism between the linear space $F_{\VV}[\emptyset;x]$ and the space of skew-symmetric elements in the free $\VV$-algebra $F_{\VV}[X]$, where $X=\{x_1,x_2,\ldots,x_n,\ldots\}$. In particular, an element $u(x)$ is an identity in the superalgebra $F_{\VV}[\emptyset;x]$ if and only if the element $({\rm Skew\,} u)$ is an identity in the free $\VV$-algebra $F_{\VV}[X]$.
\end{prop}

Let $\A=Alt[\emptyset;x]$ be the free alternative superalgebra  generated by an odd
generator $x$.   Define by induction
 \bes
 x^{[1]}=x, \ x^{[i+1]}=[x^{[i]},x]_s,\ i>0,
 \ees
and denote
$$
t=x^{[2]}, \ z^{[k]}=[x^{[k]},t],\ u^{[k]}=x^{[k]}\circ_s x^{[3]},\
k>1.
$$
The following proposition summarizes some results from \cite{ShZh-07} on the structure of $\A$.
\begin{prop}\label{prop:old} 
\begin{itemize}
\item[$(i)$]
For any   $k>0$, the element $z^{[k]}$ lies in the (super)center of $\A$; 
 
\item[$(ii)$]
Elements
\begin{gather}
t^m x^{\s},\ m+\s\geq 1,\quad
t^m (x^{[k+2]} x^{\s}), \notag \\
t^m (u^{[4k+\e]} x^{\s}),\quad t^m (z^{[4k+\e]} x^{\s}),\label{baseA}
\end{gather}
where $k>0$, $m\geq 0$;\ $\e,\,\s\in\{0,1\}$, form a base of the
superalgebra $\A$.
\end{itemize}\ctd
\end{prop}

\begin{lem}\label{lem_1}
Let $f_m$ be the skew-symmetric polynomial constructed in Section 2. For any $m>0$ there exist $\a_m,\b_m\in F$ with $\a_m\neq 0$ such that  
\bee
f_{m}= \Skew (\a_m x^{[m]}+\b_m z^{[m-2]} ).\label{id_2}
\eee
\end{lem}
\Proof The element $f_m$ is a skew-symmetric polynomial on $m$ variables in the free alternative algebra $\Alt[X]$.  Moreover, it is a Malcev element, that is, it belongs to the free special Malcev superalgebra $SMalc[\emptyset;x]$. By \cite{ShZh-07}, the superalgebra $SMalc[\emptyset;x]$ is isomorphic to the free Malcev superalgebra $Malc[\emptyset;x]$ generated by an odd element $x$, and by \cite{She-03},  the set $\{x^{[i]},\, z^{[j]}, \ i,j>0\}$ forms its base. 
By Proposition \ref{prop_u(x)}, there exists a homogeneous element $u\in SMalc[\emptyset;x]$ of degree $m$ such that $f_m=\Skew u$. 
Comparing the degrees, we have that $u=\a_m x^{[m]}+\b_m z^{[m-2]}$ for some $\a_m,\b_m\in F$. Since $\Skew z^{[m]}=0$ in any metabelian algebra, we have that $\a_m\neq 0$.
 \ctd
 
\begin{Cor}\label{cor3}
For any $m>0$ there exist $\l_m,\, \nu_m\in F$ such that $\l_m\neq 0$ and 
\bee
 \Skew x^{[m]}=\l_m f_{m}+\nu_m[\sum_{i<j}(-1)^{i+j}[f_{m-2}(x_1,\ldots,\check{x}_i,\check{x}_j,\ldots,x_{m-2}),[x_i,x_j]].\label{id_4}
\eee
\end{Cor}
\Proof
From \eqref{id_2} we obtain since $z^{[k]}\in Z(\A)$
\bes
 \Skew x^{[m]}&=&\l_m f_{m}+\nu_m(\sum_{i<j}(-1)^{i+j}[\Skew x^{[m-2]},[x_i,x_j]])\\
&=&\l_m f_{m}+\nu_m(\sum_{i<j}(-1)^{i+j}[(\l_{m-2}f_{m-2}+\nu_{m-2}\Skew z^{[m-4]},[x_i,x_j]]))\\
&=&\l_m f_{m}+\nu_m(\sum_{i<j}(-1)^{i+j}[(\l_{m-2}f_{m-2},[x_i,x_j]])),
\ees
which proves the Corollary.

\ctd

Corollaries \ref{cor2} and \ref{cor3} imply
\begin{Cor}\label{cor4}
Every $m$-generated  alternative algebra over a field of characteristic $0$ satisfies the identity  $\Skew x^{[n+3]}=0$ for $n>1+m+\left(\begin{array}{c} m\\2\end{array}\right)+\left(\begin{array}{c} m\\3\end{array}\right)$.
\end{Cor}

\begin{lem}\label{lem2}
For every natural number $m$ there exists a number $M=M(m)$ such that every alternative algebra $A$  with $m$ generators which satisfies the identity 
\bee\label{id3}
[[x,y]^2,x]=0,
\eee
satisfies also the identity
\bes
\Skew t^{4M}=0.
\ees
\end{lem}
\Proof
Let $A$ be the free $m$-generated alternative algebra in the variety defined by  identity \eqref{id3}. Since $A$ is a finitely generated $PI$-algebra, by  \cite[theorem 2.2]{She-83} its Jacobson radical $J(A)$ is nilpotent, say of degree $M$. The quotient algebra $A/J(A)$ is a subdirect sum of primitive algebras satisfying identity \eqref{id3}, which by the Kaplansky  \cite{H} and Kleinfeld \cite{ZSSS} theorems are isomorphic to matrix algebras of order no more than $2$ or to Cayley-Dickson algebras over their centers. The first algebras satisfy the identity $\Skew t^2=0$, and the second ones satisfy the identity $\Skew t^4=0$. Therefore, the quotient algebra $A/J(A)$ satisfies the identity $\Skew t^4=0$, and the algebra $A$ satisfies the identity $\Skew t^{4M}=0$.

\ctd

Now we can prove our main result.
\begin{thm*}
For every natural number $n$ there exists a number $N(n)$ such that every skew-symmetric multilinear polynomial on more then $N(n)$ variables which is zero in the free associative algebra, vanishes in any $n$-generated alternative algebra $A$.
\end{thm*}
\Proof
Let $f=f(x_1,\ldots,x_k)$ be a skew-symmetric polynomial from $Alt[X]$ which is zero in the free associative algebra $Ass[X]$. Then, as before, $f=\Skew u$, where $u=u(x)\in \A$ is a homogeneous element of degree $k$. Using base \eqref{baseA} and the fact that $f=0$ in $Ass[X]$, we conclude that $u(x)$ is a linear combination of the following elements 
\begin{enumerate}
\item $t^j(x^{[i+2]}x^{\s}),\ 2j+i+2+\s=k, i>0$
\item $t^j(u^{[4i+\e]}x^{\s}),\ 2j+4i+3+\e+\s=k,$
\item $t^j(z^{[4i+\e]}x^{\s}),\ 2j+4i+2+\e+\s=k.$
\end{enumerate}
Let $D=D(A)$ be the associator ideal of $A$, that is, the ideal generated by all associators $(a,b,c)=(ab)c-a(bc),\ a,b,c\in A$, and let $U=U(A)$ be the nucleus of $A$, that is, the maximal ideal lying in the associative center of $A$ (see \cite{ZSSS}). Then $U(A)D(A)=0$. Observe now that the second factors in products 1--3 lie in $D(A)$. Moreover, by Corollary \ref{cor4}, if $i>3+m+\left(\begin{array}{c} m\\2\end{array}\right)+\left(\begin{array}{c} m\\3\end{array}\right)$ then skew-symmetrizations of these factors are zero and $\Skew u(x)=0$.

Furthermore, by \cite{Sv}, the element $[[x,y]^2,x]\in U(A)$. Therefore, if $j\geq 4M(n)$, where $M(n)$ is the number defined in  Lemma \ref{lem2}, then $\Skew t^j\in U(A)$ and $\Skew u(x)\in U(A)D(A)=0$.  

It is clear now that we can take $N(n)=2M(n)+3+m+\left(\begin{array}{c} m\\2\end{array}\right)+\left(\begin{array}{c} m\\3\end{array}\right)$.
 
\ctd

To conclude,  we formulate the following open question:\\[1mm]
{\em Assume that an alternative algebra $A$ over a field of characteristic zero satisfies all multilinear skew-symmetric identities of a fixed degree $n$. Would $A$ satisfy for some $m$ all identities of the type
$$
\sum_{\s\in\Sigma_m} (-1)^{sgn(\s)} ((x_{\s(1)}W_1)(x_{\s(2)}W_2)\cdots (x_{\s(m)}W_m))_q=0
$$ 
for arbitrary multiplication operators $W_i$ and a fixed arrangement of brackets $q$?}

 This would be a generalization of the well known Kemer's theorem claiming that in associative algebras over a field of characteristic zero  the standard identity of degree $n$ implies the Capelly identity of a certain degree $m$ \cite{Kem}.

\end{document}